\documentclass[12pt]{article}
\usepackage{amssymb, amsmath, amsfonts, url, graphicx}

\def\0{\leqno}
\def\({\left(}
\def\){\right)}
\def\<{\left<}
\def\>{\right>}

\def\4{\subseteq }
\def\dd{\displaystyle}

\def\cale{{\cal E}}

\def\bit{\begin{itemize}}
\def\eit{\end{itemize}}
\def\barr{\begin{array}}
\def\earr{\end{array}}

\def\X#1#2{\stb{#1}{#2}{\mbox{\Huge$\times$}}}
\def\bld#1#2{{\buildrel{#1}\over{#2}}}
\def\st#1#2{{\mathrel{\mathop{#2}\limits_{#1}}{}\!}}
\def\stb#1#2#3{{\st{{#1}}{\bld{{#2}}{#3}}{}\!}}

\def\dd{\displaystyle}

\title{\bf A note on the product of element orders of finite abelian groups}
\author{Marius T\u arn\u auceanu}
\date{October 1, 2013}

\begin{document}

\maketitle

\begin{abstract}
Given a finite group $G$, we denote by $\psi\,'(G)$ the product of
element orders of $G$. Our main result proves that the restriction
of $\psi\,'$ to abelian $p$-groups of order $p^n$ is strictly
increasing with respect to a natural order on the groups relating
to the lexicographic order of the partitions of $n$. In
particular, we infer that two finite abelian groups of the same
order are isomorphic if and only if they have the same product of
element orders.
\end{abstract}

\noindent{\bf MSC (2010):} Primary 20K01; Secondary 20D60, 20D15.

\noindent{\bf Key words:} finite abelian groups, product of
element orders.

\section{Introduction}

Let $G$ be a finite group and
$$\psi(G)=\dd\sum_{x\in G}o(x),$$where $o(x)$ denotes the order of $x\in G$.
The starting point for our discussion is given by the papers
\cite{1,2} which investigate the minimum/maximum of $\psi$ on
groups of the same order. Other properties of the function $\psi$
have been studied in \cite{10} for finite abelian groups.
\smallskip

In the current note we will focus on the function
$$\psi\,'(G)=\dd\prod_{x\in G}o(x).$$In contrast with $\psi$, this is not
multiplicative, as shown by the following result.

\bigskip\noindent{\bf Proposition A.} {\it Let $G_1,G_2,...,G_k$ be finite groups having coprime orders. Then
$$\psi\,'(\X{i=1}k G_i)=\prod_{i=1}^k \psi\,'(G_i)^{n_i},$$where $n_i=\dd\prod_{^{j=1}_{j\neq i}}^k \mid\hspace{-0,5mm} G_j\hspace{-0,5mm}\mid, \hspace{1mm}i=1,2,...,k$. In particular,
if $G$ is a finite nilpotent group of order
$n=p_1^{\alpha_1}p_2^{\alpha_2}\cdots p_k^{\alpha_k}$,
$G_1,G_2,...,G_k$ are the Sylow subgroups of $G$ and
$n_i=n/p_i^{\alpha_i}, \hspace{1mm}i=1,2,...,k$, then
$$\psi\,'(G)=\prod_{i=1}^k \psi\,'(G_i)^{n_i}.$$}

By Proposition A we infer that the computation of $\psi\,'(G)$ for
nilpotent groups is reduced to $p$-groups and explicit formulas
can be given in several particular cases. One of them consists of
abelian groups, for which Corollary 4.4 of \cite{9} leads to the
following theorem.

\smallskip\noindent{\bf Theorem B.} {\it Let $G=\X{i=1}k\,
\mathbb{Z}_{p^{\alpha_i}}$ be a finite abelian $p$-group, where
$1\leq\alpha_1\leq\alpha_2\leq...\leq\alpha_k$. Then
$$\psi\,'(G)=p^{\hspace{2mm}\alpha_k p^{\alpha_1+\alpha_2+...+\alpha_k}-\sum_{i=0}^{\alpha_k\hspace{-0,5mm}-1}p^i f_{(\alpha_1,\alpha_2,...,\alpha_k)}(i)},\0(1)$$where
$$f_{(\alpha_1,\alpha_2,...,\alpha_k)}(i)=\left\{\barr{lll}
p^{(k-1)i},&\mbox{if }&0\le i\le\alpha_1\vspace*{1,5mm}\\
p^{(k-2)i+\alpha_1},&\mbox{if }&\alpha_1\le i\le\alpha_2\\
p^{(k-3)i+\alpha_1+\alpha_2},&\mbox{if }&\alpha_2\le i\le\alpha_3\\
\vdots\\
p^{\alpha_1+\alpha_2+...+\alpha_{k-1}},&\mbox{if
}&i\geq\alpha_{k-1}\,(where\hspace{1mm} \alpha_{k-1}=1
\hspace{1mm}if\hspace{1mm}k=1).\earr\right.$$}
\smallskip

We exemplify (1) by computing $\psi\,'(G)$ for cyclic $p$-groups
and for rank two abelian $p$-groups.

\bigskip\noindent{\bf Example.} We have:
\begin{itemize}
\item[1.] $\psi\,'(\mathbb{Z}_{p^{\alpha}})=p^{\,\frac{\alpha p^{\alpha+1}-(\alpha+1)p^{\alpha}+1}{p-1}}$\,;
\item[2.] $\psi\,'(\mathbb{Z}_{p^{\alpha}}\times\mathbb{Z}_{p^{\beta}})=p^{\,\frac{\beta p^{\alpha+\beta+2}-p^{\alpha+\beta+1}-(\beta+1)p^{\alpha+\beta}+p^{2\alpha+1}+1}{p^2-1}}\,.$
\end{itemize}
\smallskip

Given a positive integer $n$, it is well-known that there is a
bijection between the set of types of abelian groups of order
$p^n$ and the set $P_n=\{(x_1,x_2,...,x_n)\in\mathbb{N}^n\mid
x_1\geq x_2\geq...\geq x_n,\, x_1+x_2+...+x_n=n\}$ of partitions
of $n$. Namely,
$$\X{i=1}k\,\mathbb{Z}_{p^{\alpha_i}}(\mbox{with }1\leq\alpha_1\leq\alpha_2\leq...\leq\alpha_k\mbox{ and}\sum_{i=1}^k \alpha_i=n)\mapsto (\alpha_k,...,\alpha_1,\hspace{-3mm}\underbrace{0,...,0}_{n-k\,
\rm positions }\hspace{-2mm})\,.$$Moreover, recall that $P_n$ is
totally ordered under the lexicographic order $\preceq$\,, where
$$(x_1,x_2,...,x_n)\hspace{-1mm}\prec\hspace{-1mm}(y_1,y_2,...,y_n)\hspace{-1mm}\Longleftrightarrow\hspace{-1mm} \left\{\barr{lll}
x_1=y_1,...,x_m=y_m\\
\mbox{and}\\
x_{m+1}{<} y_{m+1} \mbox{ for some }
m\hspace{-1mm}\in\hspace{-1mm}\left\{0,1,...,n{-}1\right\}.\earr\right.$$This
induces a total order on the set of types of abelian $p$-groups of
order $p^n$. Our next theorem shows that the restriction of
$\psi\,'$ on this set is strictly increasing.

\bigskip\noindent{\bf Theorem C.} {\it Let $G=\X{i=1}k\,\mathbb{Z}_{p^{\alpha_i}}$ and $H=\X{j=1}r\,\mathbb{Z}_{p^{\beta_j}}$ be two finite abelian
$p$-groups of order $p^n$. Then
$$\psi\,'(G)<\psi\,'(H)\Longleftrightarrow(\alpha_k,...,\alpha_1,\hspace{-2mm}\underbrace{0,...,0}_{n-k\,
\rm positions
}\hspace{-2mm})\prec(\beta_r,...,\beta_1,\hspace{-2mm}\underbrace{0,...,0}_{n-r\,
\rm positions }\hspace{-2mm})\,.\0(2)$$}

Since a strictly increasing function is injective, by Theorem C we
infer the following corollary.

\bigskip\noindent{\bf Corollary D.} {\it Two finite abelian $p$-groups of order $p^n$ are isomorphic if and only if
they have the same product of element orders.}
\bigskip

This can be extended to arbitrary finite abelian groups, according
to Proposition A.

\bigskip\noindent{\bf Theorem E.} {\it Two finite abelian groups of the same order are isomorphic if and only if
they have the same product of element orders.}

\bigskip\noindent{\bf Remark.} The above two results are not true
for \textit{arbitrary} finite abelian groups. For example, we have
$\psi\,'(\mathbb{Z}_4\times\mathbb{Z}_3^2)=\psi\,'(\mathbb{Z}_2^4\times\mathbb{Z}_3)$,
but obviously the groups $\mathbb{Z}_4\times\mathbb{Z}_3^2$ and
$\mathbb{Z}_2^4\times\mathbb{Z}_3$ are not isomorphic.
\bigskip

Finally, we associate to a finite (abelian) group
$G=\{x_1,x_2,...,x_n\}$ the polynomial
$$P_G=\dd\prod_{i=1}^n\,(X-o(x_i))\in \mathbb{Z}[X].$$Recall that
if $G$ and $H$ are two finite abelian groups for which $P_G=P_H$
(that is, $G$ and $H$ have the same element orders with the same
multiplicities), then $G\cong H$ by Theorem 5 of \cite{8}. In
order to improve this result, we construct the quantities
$$\psi_k(G)=\sum_{1\leq i_1<i_2<...<i_k\leq n}o(x_{i_1})o(x_{i_2})\cdots
o(x_{i_k}), \hspace{1mm}k=1,2,...,n.$$Obviously, explicit formulas
for $\psi_k(G)$, $k=1,2,...,n$, can be given by using Corollary
4.4 of \cite{9}. We also observe that $\psi_1(G)=\psi(G)$ and
$\psi_n(G)=\psi\,'(G)$. Inspired by Conjecture 6 of \cite{10} and
the above Theorem E, we came up with the following conjecture,
which we have verified by GAP for many values of $k$ and $n$.

\bigskip\noindent{\bf Conjecture F.} {\it Let $G$ and $H$ be two finite abelian groups of order $n$. Then for every $k\in\{1,2,...,n\}$, we have
$G\cong H$ if and only if $\psi_k(G)=\psi_k(H)$.}
\bigskip

Most of our notation is standard and will not be repeated here.
For basic notions and results of group theory we refer the reader
to \cite{4,5}. Other interesting papers on the above topic are
\cite{3,6,7}.

\section{Proof of Theorem C}

We first remark that it suffices to verify (2) only for
consecutive partitions of $n$, because $P_n$ is totally ordered.
\smallskip

Assume that
$(\alpha_k,...,\alpha_1,0,...,0)\prec(\beta_r,...,\beta_1,0,...,0)$
and let $s\in\{1,2,...,\newline r-1\}$ such that
$\beta_1=\beta_2=\cdots=\beta_s<\beta_{s+1}$. We distinguish the
following two cases.

\medskip\noindent{\bf \hspace{20mm}Case 1.} $\beta_1\geq 2$

\noindent Then $(\alpha_k,...,\alpha_1,0,...,0)$ is of type
$(\beta_r,...,\beta_2,\beta_1-1,1,0,...,0)$, i.e. $k=r+1$,
$\alpha_1=1$, $\alpha_2=\beta_1-1$ and $\alpha_i=\beta_{i-1}$ for
$i=3,4,...,r+1$. It is easy to see that
$f_{(\alpha_1,\alpha_2,...,\alpha_k)}(i)=f_{(\beta_1,\beta_2,...,\beta_r)}(i)
\mbox{ for all} \hspace{1mm}i\geq\beta_1$. One obtains
$$\psi\,'(G)<\psi\,'(H)\Longleftrightarrow \dd\sum_{i=1}^{\beta_r-1} p^i f_{(\beta_1,\beta_2,...,\beta_r)}(i) <\dd\sum_{i=1}^{\alpha_k-1}p^i f_{(\alpha_1,\alpha_2,...,\alpha_k)}(i)$$
$$\hspace{29mm}\Longleftrightarrow\dd\sum_{i=1}^{\beta_1-1} p^i
f_{(\beta_1,\beta_2,...,\beta_r)}(i) <\dd\sum_{i=1}^{\beta_1-1}p^i
f_{(\alpha_1,\alpha_2,...,\alpha_k)}(i)$$and the last inequality
is true because
$$f_{(\beta_1,\beta_2,...,\beta_r)}(i){=}p^{(r-1)i}{<}p^{(r-1)i+1}{=}p^{(k-2)i+1}{=}f_{(\alpha_1,\alpha_2,...,\alpha_k)}(i), i=1,2,...,\beta_1{-}1\,.$$

\medskip\noindent{\bf \hspace{20mm}Case 2.} $\beta_1=1$

\noindent Then $(\alpha_k,...,\alpha_1,0,...,0)$ is of type
$(\beta_r,...,\beta_{s+1}-1,\beta'_t,\beta'_{t-1},...,\beta'_1,0,...,0)$,
where
$\beta_{s+1}-1\geq\beta'_t\geq\beta'_{t-1}\geq...\geq\beta'_1\geq
1$ and $\beta'_t+\beta'_{t-1}+...+\beta'_1=s+1$. We infer that
$f_{(\alpha_1,\alpha_2,...,\alpha_k)}(i)=f_{(\beta_1,\beta_2,...,\beta_r)}(i)
\mbox{ for all} \hspace{1mm}i\geq\beta_{s+1}$. So, we can suppose
that $s=r-1$, i.e.
$$(\alpha_k,...,\alpha_1,0,...,0)=(\beta_{r}-1,\beta'_t,\beta'_{t-1},...,\beta'_1,0,...,0).$$One
obtains
$$\psi\,'(G){<}\psi\,'(H){\Longleftrightarrow} (\beta_r{-}1)p^n{-}\hspace{-2mm}\dd\sum_{i=0}^{\beta_r-2}p^i f_{(\alpha_1,\alpha_2,...,\alpha_k)}(i){<}\beta_r p^n{-}\hspace{-2mm}\dd\sum_{i=0}^{\beta_r-1}p^i f_{(\beta_1,\beta_2,...,\beta_r)}(i)$$
$$\hspace{23mm}{\Longleftrightarrow}\, p^n{-}\hspace{-1mm}\dd\sum_{i=0}^{\beta_r-1}p^i f_{(\beta_1,\beta_2,...,\beta_r)}(i)+\hspace{-1mm}\dd\sum_{i=0}^{\beta_r-2}p^i
f_{(\alpha_1,\alpha_2,...,\alpha_k)}(i)>0\,.\hspace{3mm}(*)$$Since
$$f_{(\beta_1,\beta_2,...,\beta_r)}(i)=\left\{\barr{lll}
p^{(r-1)i},&\mbox{ if }&0\le i\le 1\vspace*{1,5mm}\\
p^{r-1},&\mbox{ if }&i\geq 1\,,\earr\right.$$we easily get
$$\dd\sum_{i=0}^{\beta_r-1}p^i
f_{(\beta_1,\beta_2,...,\beta_r)}(i)=1+p^{r-1}\dd\frac{p^{\beta_r}-p}{p-1}=1+\dd\frac{p^n-p^r}{p-1}<p^n$$and
therefore $(*)$ is true. \smallskip

Conversely, assume that
$(\alpha_k,\alpha_{k-1},...,\alpha_1,0,...,0){\succeq}(\beta_r,\beta_{r-1},...,\beta_1,0,...,0)$.
If these partitions are equal, then $G\cong H$, so
$\psi\,'(G)=\psi\,'(H)$. Otherwise,
$(\alpha_k,\alpha_{k-1},...,\alpha_1,0,...,0){\succ}(\beta_r,\beta_{r-1},...,\beta_1,0,...,0)$,
and the first part of the proof implies $\psi\,'(G)>\psi\,'(H)$,
as desired.
$\scriptstyle\Box$
\bigskip

\bigskip\noindent{\bf Acknowledgements.} The author is grateful to the reviewers for
their remarks which improve the previous version of the paper.

\vspace*{5ex}\small

\hfill
\begin{minipage}[t]{5cm}
Marius T\u arn\u auceanu \\
Faculty of  Mathematics \\
``Al.I. Cuza'' University \\
Ia\c si, Romania \\
e-mail: {\tt tarnauc@uaic.ro}
\end{minipage}

\end{document}